\documentclass[12pt, twoside, leqno]{article}
\usepackage{hyperref}

\usepackage{amsmath,amsthm}
\usepackage{amssymb}

\usepackage{enumitem}

\usepackage{graphicx}

\usepackage[T1]{fontenc}

\usepackage{amsfonts}
\usepackage{float}
\usepackage{cancel}

\theoremstyle{definition}

\numberwithin{equation}{section}

\begin{document}


\baselineskip=17pt


\title{THE BEST CONSTANT FOR INEQUALITY INVOLVING SUM OF THE RECIPROCALS AND PRODUCT OF POSITIVE NUMBERS WITH UNIT SUM}

\author{Yagub N. Aliyev\\
School of IT and Engineering\\ 
ADA University\\
Ahmadbey Aghaoglu str. 61 \\
Baku 1008, Azerbaijan\\
E-mail: yaliyev@ada.edu.az}
\date{}

\maketitle


\renewcommand{\thefootnote}{}



\renewcommand{\thefootnote}{\arabic{footnote}}
\setcounter{footnote}{0}


\begin{abstract}
In the paper we study a special parameter containing algebraic inequality involving sum of reciprocals and product of positive real numbers whose
sum is 1. We determine the best values of the parameter using a new optimization argument. In the case of three numbers the algebraic inequality have some interesting geometric applications involving a
generalization of Euler's inequality about the ratio of radii of circumscribed
and inscribed circles of a triangle.
\end{abstract}

\section{Introduction}
Consider the inequality
$$
\sum\limits_{i = 1}^n{\frac{1}{x_{i}}}\geq \frac{\lambda
}{1+n^{n-2}(\lambda -n^{2})\prod\limits_{i=1}^n{x_{i}}}, \eqno(1)
$$
where $x_{1},x_{2},\ldots ,x_{n}>0; \sum\limits_{i =
1}^n{x_{i}}=1$, for $n\geq 2$. Here $\lambda >0$ is a real number and it is asked to find the best (maximal possible)  $\lambda $ for each  $n$. If such $\lambda $ exists then we will denote it by $\lambda_n$. 
Note that the right hand side of the inequality (1)
$$
f(\lambda)= \frac{\lambda
}{1+n^{n-2}(\lambda -n^{2})\prod\limits_{i=1}^n{x_{i}}},
$$
where $x_{1},x_{2},\ldots ,x_{n}>0; \sum\limits_{i =
1}^n{x_{i}}=1$, is a non-decreasing function of $\lambda >0$. So if (1) is true for a certain $\lambda =\lambda _n$ then it is also true for all $0<\lambda \leq \lambda _n$.

By Cauchy-Schwarz inequality $\sum\limits_{i = 1}^n{\frac{1}{x_{i}}}\geq n^2 =f(n^2)$. Since the inequality holds true for $\lambda = n^2$, it also holds true for all $0<\lambda \leq n^2$. So, the best constant $\lambda =\lambda _n$, if exists, satisfies $\lambda_n \geq n^2$.

\textit{Case $n=2$.} For the case $n=2$ there is no best constant. If $n=2$  then we obtain the inequality
$$
\frac{1}{x_{1}}+\frac{1}{x_{2}}\geq \frac{\lambda
}{1+(\lambda -2^{2})x_{1}x_{2}},
$$
where $x_{1},x_{2}>0; x_{1}+x_{2}=1$. This inequality is true for any $\lambda >0$. Indeed, if we multiply both sides by $(1+(\lambda -4)x_{1}x_{2})$, then we obtain
$$
\frac{1}{x_{1}}+\frac{1}{x_{2}}+(\lambda -4)(x_{1}+x_{2})\geq \lambda.
$$
Since $x_{1}+x_{2}=1$, the parameter $\lambda $ cancels out, and we obtain
$$
\frac{1}{x_{1}}+\frac{1}{x_{2}}\geq 4,
$$
which is always true.

\textit{Case $n=3$.} For case $n=3$ the best constant is $\lambda_3=25$. We obtain the inequality
$$
\frac{1}{x_1}+\frac{1}{x_2}+\frac{1}{x_3} \geq \frac{\lambda
}{1+3(\lambda -9)x_1x_2x_3},
$$
where $x_{1},x_{2},x_{3}>0; x_{1}+x_{2}+x_{3}=1$. This inequality is true only for $0<\lambda \leq 25$.
We can show this by substituting 
$x_1=x_2=\frac{1}{4}$, $x_3=\frac{1}{2}$ in this inequality. On the other hand we can prove that
$$
 \frac{1}{x_1}+\frac{1}{x_2}+\frac{1}{x_3}\geq \frac{25}{1+48x_1x_2x_3},
$$
holds true. So, $\lambda =25$ is the maximum possible value for this inequality (see \cite{aliyev1}). In the solution to problem \cite{aliyev1} it was noticed by David B. Leep that the case $\lambda =25$ is equivalent to a more general inequality $s_1^3s_2+48s_2s_3-25s_1^2s_3\ge 0$ for symmetric polynomials $s_1=x_1+x_2+x_3$, $s_2=x_1x_2+x_2x_3+x_3x_1$, $s_3=x_1x_2x_3$, which can also be writen as
$$
x_1(x_2-x_3)^2(3x_1-x_2-x_3)^2
$$
$$
+x_2(x_1-x_3)^2(3x_2-x_1-x_3)^2+x_3(x_2-x_1)^2(3x_3-x_2-x_1)^2\ge 0,
$$
making case $n=3$ almost trivial. Inequality (1) can also be written using symmetric polynomials but as the results for cases $n=4$ and $n=5$ below suggest there is no simple solution for $n>3$. Let $$s_1=\sum_{i=1}^{n}{x_i},\ s_{n-1}=\sum_{i=1}^{n}{\prod_{j=1, j\ne i}^{n}{x_j}},\ s_n=\prod_{i=1}^{n}{x_i}.$$
If $\lambda >0$, then (1) is equivalent to inequality
$$
s_1^ns_{n-1}+n^{n-2}(\lambda -n^{2})s_{n-1}s_n-\lambda s_1^{n-1}s_n\ge 0,
$$
homogeneous with respect to variables $x_1,\ldots, x_n$.

There are some geometric applications of case $n=3$. The inequality
$$
\frac{R}{r}\geq 2+\mu
\frac{(a-b)^{2}+(b-c)^{2}+(c-a)^{2}}{(a+b+c)^{2}}, \eqno(2)
$$
where $R$ and $r$ are respectively the circumradius and
inradius, and $a,b,c$ are sides of a triangle, holds true if $\mu
\leq 8$. Indeed, substituting $a=b=3$, $c=2$ and the corresponding values of
$R=\frac{9}{4\sqrt{2}}$ and $r=\frac{\sqrt{2}}{2}$ in (2) we obtain $\mu \leq 8$. So, again, if we will prove (2) for $\mu =8$ then $\mu =8$ will the best constant for the inequality (2). For $\mu =8$ we obtain
$$
\frac{R}{r}\geq 2+8
\frac{(a-b)^{2}+(b-c)^{2}+(c-a)^{2}}{(a+b+c)^{2}},
$$
which is a refinement of Euler's inequality $\frac{R}{r}\geq 2$ and follows directly from the case $n=3$ (see \cite{aliyev2},  \cite{aliyev3}).

Another geometric application is the following inequality about the sides $a,b,c$ of a triangle which follows directly from the case $n=3$ (see \cite{aliyev4}):
$$
\frac{a^3}{b+c-a}+\frac{b^3}{a+c-b}+\frac{c^3}{a+b-c}+7(ab+bc+ca) \ge 8(a^2+b^2+c^2).
$$
This inequality can also be written as a quintic inequality of symmetric polynomials
$$
9\sum^{3}{a^5} - 15\sum^{6}{a^4b}+ 6\sum^{6}{a^3b^2} + 25\sum^{3}{a^3bc} - 16\sum^{3}{ab^2c^2} \ge 0,
$$
which is a special case ($v=3$) of the following inequality mentioned in \cite{harris} (see p. 244, where $v=u+1$)
$$
v^2\sum^{3}{a^5} - v(v+2)\sum^{6}{a^4b}+ 2v\sum^{6}{a^3b^2}
$$
$$
+ (v+2)^2\sum^{3}{a^3bc} - 4(v+1)\sum^{3}{ab^2c^2} \ge 0.
$$
This general inequality is also easily proved if we put $a=x_2+x_3$, $b=x_1+x_3$, $c=x_1+x_2$, and simplfy to obtain
$$
4x_1(x_2-x_3)^2(vx_1-x_2-x_3)^2
$$
$$
+4x_2(x_1-x_3)^2(vx_2-x_1-x_3)^2+4x_3(x_2-x_1)^2(vx_3-x_2-x_1)^2\ge 0.
$$
Similar quartic and sextic inequalities were studied in series of works by J.F. Rigby, O. Bottema and J.T. Groenman (see \cite{bottema}, \cite{rigby} and their references).

One more geometric application of case $n=3$ is about the areas of triangles and needs introduction of some notations. Let $M$ be a point in a triangle $ABC$. Extend lines $AM$, $BM$, and $CM$ to intersect the sides of triangle $ABC$ at $A_0$, $B_0$, and $C_0$, respectively (see Fig. 1). Next construct the parallel to $A_0C_0$ through $M$, which intersects $BA$ and $BC$ at $C_1$ and $A_2$, respectively. Analogously, draw the parallel through $M$ to $B_0A_0$ (and to $B_0C_0$) to find $A_1$ and $B_2$ (and $B_1$ and $C_2$). Denote
$$
T_1=[MC_1B_2], T_2=[MA_1C_2], T_3=[MB_1A_2],
$$
$$
S_1=[MA_1A_2], S_2=[MB_1B_2], S_3=[MC_1C_2],
$$
$$
P_1=[AB_2C_1], P_2=[BC_2A_1], P_3=[CA_2B_1],
$$
where the square brackets stand for the area of the triangles (see \cite{aliyev4}, \cite{aliyev5}). Then $$P_1+P_2+P_3+7(S_1+S_2+S_3)\ge 8(T_1+T_2+T_3).$$ 

\begin{figure}[htbp]
\centerline{\includegraphics[scale=.5]{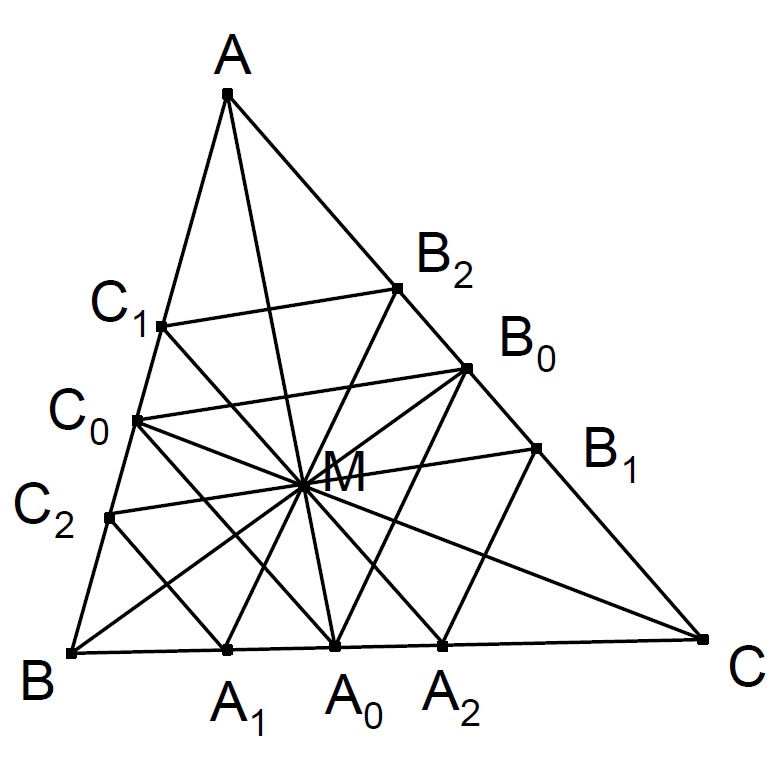}\includegraphics[scale=0.8]{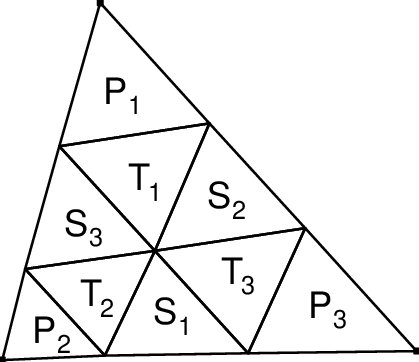}}
\label{fig1}
\caption{Geometric application of case $n=3$.}
\end{figure}

\textit{Case $n=4$.} For the case $n=4$ the best constant is $\lambda_4=\frac{582\sqrt{97}-2054}{121}\approx 30.4$ (see \cite{aliyev2}). In this case we obtain
$$
\frac{1}{x_1}+\frac{1}{x_2}+\frac{1}{x_3}+\frac{1}{x_4} \geq \frac{\lambda }{1+16(\lambda -16)x_1x_2x_3x_4},
$$
where $x_{1},x_{2},x_{3},x_{4}>0; x_{1}+x_{2}+x_{3}+x_{4}=1$. Again this inequality is true only for $\lambda \leq \frac{582\sqrt{97}-2054}{121}$. Indeed, substituting in this inequality
$x_1=x_2=x_3=\frac{5+\sqrt{97}}{72}$, $x_4=\frac{19-\sqrt{97}}{24}$, we obtain $0<\lambda \leq \frac{582\sqrt{97}-2054}{121}$. On the other hand we can prove that the inequality holds true for $\lambda = \frac{582\sqrt{97}-2054}{121}$.
So, $\lambda = \frac{582\sqrt{97}-2054}{121}$ is the maximum possible value for this inequality.

\textit{Case $n=5$.} For the case $n=5$ we obtain the inequality
$$
\frac{1}{x_1}+\frac{1}{x_2}+\frac{1}{x_3}+\frac{1}{x_4}+\frac{1}{x_5} \geq \frac{\lambda }{1+125(\lambda -25)x_1x_2x_3x_4x_5},
$$
where $x_{1},x_{2},x_{3},x_{4},x_{5}>0; x_{1}+x_{2}+x_{3}+x_{4}+x_{5}=1$, and it was conjectured in \cite{aliyev2} that  the best constant is
$$
\lambda _5=\frac{12933567-93093\sqrt{22535}}{4135801}\alpha
$$
$$
+\frac{17887113+560211\sqrt{22535}}{996728041}\alpha
^{2}-\frac{288017}{17161}\approx 40.09,
$$
where $\alpha =\sqrt[3]{8119+48\sqrt{22535}}$. This conjecture for $\lambda _5$ will be proved in the current paper. Also it will be proved that the equality cases in this inequality occur when $x_1=x_2=x_3=x_4=x_5= \frac{1}{5}$ and when, for example,
$x_1=x_2=x_3=x_4=x=\frac{\alpha}{240} + \frac{241}{240\alpha} + \frac{7}{240}\approx 0.173$, $x_5=1-4x\approx 0.308$.

\textit{Case $n=6$.} This case was not studied before. Using Maple the exact value of $\lambda _6$ is calculated.
Case $n=6$ is possibly the last case for which these calculations of the exact value are possible.

\textit{Case $n\ge 7$.} In view of the fact that quintic and higher order equations are in general not solvable in radicals, it is unlikely that there is a precise formula for the best constant in the cases $n\ge 7$. Therefore, for the greater values of $n$ ($n\ge 7$), instead of the exact value, it is reasonable to find some bounds or approximations for $\lambda_n$. In the current paper it is proved that
$$
\frac{n^3}{n-1}\leq \lambda _n\leq  \frac{n^3}{n-2}. \eqno(3)
$$
Some possible improvements for this symmetric double inequality are also discussed.

It is interesting to compare the results of the current paper with the results for the following similar inequality
$$\sum\limits_{i =
1}^n {\frac{1}{x_{i}}}\leq \nu +\frac{n^{2}-\nu
}{n^{n}\prod \limits_{i = 1}^n {x_{i}}}, \eqno(4)
$$
where $x_{1},x_{2},\ldots ,x_{i}>0;$ $\sum\limits_{i = 1}^n
{x_{i}}=1.$ The best constant $\nu_n$ for this inequality is known for all $n>1$. See Corollary 2.13 in \cite{mitev} where it is proved that $\nu\le\nu_n
=n^{2}-\frac{n^{n}}{(n-1)^{n-1}}$.
In particular, if $\nu =0$ then we obtain
$$\sum\limits_{i =
1}^n {\frac{1}{x_{i}}}\leq \frac{1
}{n^{n-2}\prod \limits_{i = 1}^n {x_{i}}}, \eqno(5)
$$
with equality case possible only when $x_1=\ldots=x_n=\frac{1}{n}$. Inequality (5) also follows from the following inequality for $E_i=\frac{1}{{\binom{n}{i}}}s_i$ (averages of $s_i$),
$$
E_1^{\alpha_1}\ldots E_n^{\alpha_n}\le E_1^{\beta_1}\ldots E_n^{\beta_n},
$$
which holds if and only if
$$
\alpha_m + 2\alpha_{m+1}+\ldots+(n- m + 1)\alpha_n\ge \beta_m + 2\beta_{m+1}+\ldots+(n- m + 1)\beta_n,
$$
for each $1\le m\le n$ (see \cite{mitev}, Theorem 1.1; \cite{hardy}, p. 94, item 77). Indeed, it is sufficient to note that inequality (5) can be written as $E_{n-1}\le E_1^{n-1}$. This means that the above conditions for $\alpha_i,\beta_i$ ($i=1,\ldots,n$) are satisfied as
$$
\alpha_1=\ldots=\alpha_{n-2}=0,\ \alpha_{n-1}=n-1,\ \alpha_{n}=0,\ \beta_1=n-1,\ \beta_2=\ldots=\beta_{n}=0.
$$
Since (5) will be essential in the following text, an independent proof of (5) and some generalizations will be given in Appendix. Note also that 
$$
\lim_{\lambda\rightarrow +\infty}\frac{\lambda
}{1+n^{n-2}(\lambda -n^{2})\prod\limits_{i=1}^n{x_{i}}}=\frac{1
}{n^{n-2}\prod \limits_{i = 1}^n {x_{i}}}.
$$
Using this and by comparing (1) and (5), we obtain that if $n>2$, then $\lambda_n <+\infty$.

Special cases $n=3$ and $n=4$ of inequality (4) are also of interest for comparison with the corresponding cases of inequality (1). If $n=3$ then the
best constant inequality is $\frac{1}{x_1}+\frac{1}{x_2}+\frac{1}{x_3}\leq
\frac{9}{4}+\frac{1}{4x_1x_2x_3},$
where $x_1,x_2,x_3>0;$ $x_1+x_2+x_3=1$. Surprisingly, this inequality is also equivalent to a
geometric inequality. One can show that it simplifies to $p^{2}\geq 16Rr-5r^{2},$
where $p$ is the semiperimeter of a triangle. The last geometric inequality also
follows from the formula for the distance between the
incenter $I$ and the centroid $G$ of a triangle: $\left| IG\right| ^{2}=\frac{1}{9}(p^{2}+5r^{2}-16Rr)$ (see \cite{aliyev2}). 
If $n=4$ then the
best constant inequality is $\frac{1}{x_1}+\frac{1}{x_2}+\frac{1}{x_3}+\frac{1}{x_4}\leq
\frac{176}{27}+\frac{1}{27x_1x_2x_3x_4},$
where $x_1,x_2,x_3,x_4>0;$ $x_1+x_2+x_3+x_4=1$ (see \cite{sato}, Example 3).

The literature about symmetric polynomial inequalities is extensive \cite{mitrin1}, \cite{mitrin2}, \cite{mitrin3}, \cite{nicul1}, \cite{nicul2},  \cite{rosset},  \cite{tao}. Some of the results of the current paper were presented at Maple Conference 2021 \cite{aliyev7}.

\section{Main Results}

Let us consider all cases for $n\ge 3$ in a unified way. Assume first that $\left(x_{1},x_{2},\ldots ,x_{n}\right)\ne \left(\frac{1}{n},\frac{1}{n},\ldots,\frac{1}{n}\right)$. Then by using (5), inequality (1) can be written as
$$
\frac{n^2\left(1-n^n\prod\limits_{i=1}^n{x_{i}}\right)}{\frac{n^2}{\sum\limits_{i = 1}^n{\frac{1}{x_{i}}}}-n^n\prod\limits_{i=1}^n{x_{i}}}\geq \lambda , \eqno(6)
$$
where $x_{1},x_{2},\ldots ,x_{n}>0; \sum\limits_{i =
1}^n{x_{i}}=1$, for $n\geq 3$. Let us denote the left side of (6) by $g(x_1,\ldots,x_n)$, which is defined for all points of bounded set
$$
C=\{\textbf{x}|\textbf{x}=(x_{1},x_{2},\ldots ,x_{n});\ x_{1},x_{2},\ldots ,x_{n}\ge 0;\sum\limits_{i =
1}^n{x_{i}}=1\},
$$
except point $P_0\left(\frac{1}{n},\frac{1}{n},\ldots,\frac{1}{n}\right)$. For the points of boundary
$$
\partial C=\{\textbf{x}|\textbf{x}=(x_{1},x_{2},\ldots ,x_{n});\ x_{1},x_{2},\ldots ,x_{n}\ge 0;\sum\limits_{i =
1}^n{x_{i}}=1,\ \prod\limits_{i=1}^n{x_{i}}=0\},
$$
function $g$ is undefined and, obviously, for each $i=1,\ldots,n$,
$$
\lim_{x_i\rightarrow 0} g(x_1,\ldots,x_n) =+\infty.
$$
\textbf{Lemma 1.} If $x>0$, then
$$
\lim_{(x_1,x_2,\ldots,x_n)\rightarrow (x,x,\ldots,x)} \frac{\left(\sum\limits_{i =
1}^n{x_{i}}\right)-\frac{n^2}{\sum\limits_{i = 1}^n{\frac{1}{x_{i}}}}}{\left(\sum\limits_{i =
1}^n{x_{i}}\right)^n-n^n\prod\limits_{i=1}^n{x_{i}}}=\frac{2}{n^nx^{n-1}}.
$$
\textit{Proof.} The limit can be interpreted as a single variable limit if we take $$(x_1,x_2,\ldots,x_n)=(x+\gamma_1t,x+\gamma_2t,\ldots,x+\gamma_nt),$$ where not all constants $\gamma_i$ are equal and $t\rightarrow 0$. So, we calculate
$$
\lim_{t\rightarrow 0} \frac{\left(\sum\limits_{i =
1}^n{\left(x+\gamma_it\right)}\right)-\frac{n^2}{\sum\limits_{i = 1}^n{\frac{1}{x+\gamma_it}}}}{\left(\sum\limits_{i =
1}^n{\left(x+\gamma_it\right)}\right)^n-n^n\prod\limits_{i=1}^n{\left(x+\gamma_it\right)}}
$$
$$
=\lim_{t\rightarrow 0} \frac{\left(\sum\limits_{i =
1}^n{\gamma_i}\right)-\frac{n^2}{\left(\sum\limits_{i = 1}^n{\frac{1}{x+\gamma_it}}\right)^2}\sum\limits_{i = 1}^n{\frac{\gamma_i}{(x+\gamma_it)^2}}}{n\left(\sum\limits_{i =
1}^n{\left(x+\gamma_it\right)}\right)^{n-1}\left(\sum\limits_{i =
1}^n{\gamma_i}\right)-n^n\left({\prod\limits_{i=1}^n{\left(x+\gamma_it\right)}}\right)\left(\sum\limits_{i =
1}^n{\frac{\gamma_i}{x+\gamma_it}}\right)}
$$
$$
=\lim_{t\rightarrow 0} \frac{\frac{-2n^2}{\left(\sum\limits_{i = 1}^n{\frac{1}{x+\gamma_it}}\right)^3}\left(\sum\limits_{i = 1}^n{\frac{\gamma_i}{(x+\gamma_it)^2}}\right)^2+\frac{n^2}{\left(\sum\limits_{i = 1}^n{\frac{1}{x+\gamma_it}}\right)^2}\left(\sum\limits_{i = 1}^n{\frac{2\gamma_i^2}{(x+\gamma_it)^3}}\right)
}{n(n-1)\left(\sum\limits_{i =
1}^n{\left(x+\gamma_it\right)}\right)^{n-2}\left(\sum\limits_{i =
1}^n{\gamma_i}\right)^2-n^n\left({\prod\limits_{i=1}^n{\left(x+\gamma_it\right)}}\right)\left(
\left(\sum\limits_{i =
1}^n{\frac{\gamma_i}{x+\gamma_it}}\right)^2-\sum\limits_{i =
1}^n{\frac{\gamma_i^2}{(x+\gamma_it)^2}}\right)}
$$
$$
=\frac{\frac{-2n^2}{\left( \frac{n}{x}\right)^3}\left(\frac{\sum\limits_{i = 1}^n{\gamma_i}}{x^2}\right)^2+\frac{n^2}{\left(\frac{n}{x}\right)^2}\frac{\sum\limits_{i = 1}^n{2\gamma_i^2}}{x^3}
}{n(\cancel{n}-1)\left(nx\right)^{n-2}\left(\sum\limits_{i =
1}^n{\gamma_i}\right)^2-n^nx^n\left(
\cancel{\left(\frac{\sum\limits_{i = 1}^n{\gamma_i}}{x}\right)^2}-\frac{\sum\limits_{i = 1}^n{\gamma_i^2}}{x^2}\right)}=\frac{2}{n^nx^{n-1}},
$$
where we used L'Hôpital's rule twice and the fact that $n\sum\limits_{i = 1}^n{\gamma_i^2}>\left(\sum\limits_{i =
1}^n{\gamma_i}\right)^2$ (Cauchy-Schwarz inequality, the equality case is not possible as not all $\gamma_i$ are equal). Proof is complete.

In particular, if $\sum\limits_{i =
1}^n{x_{i}}=1$ then $x=\frac{1}{n}$, and therefore, by Lemma 1,
$$
\lim_{(x_1,x_2,\ldots,x_n)\rightarrow \left(\frac{1}{n},\frac{1}{n},\ldots,\frac{1}{n}\right)
} \frac{n^2\left(1-n^n\prod\limits_{i=1}^n{x_{i}}\right)}{\frac{n^2}{\sum\limits_{i = 1}^n{\frac{1}{x_{i}}}}-n^n\prod\limits_{i=1}^n{x_{i}}}
$$
$$
=\frac{n^2}{1-\lim_{(x_1,x_2,\ldots,x_n)\rightarrow \left(\frac{1}{n},\frac{1}{n},\ldots,\frac{1}{n}\right)}{\frac{1-\frac{n^2}{\sum\limits_{i = 1}^n{\frac{1}{x_{i}}}}}{1-n^n\prod\limits_{i=1}^n{x_{i}}}}}=\frac{n^2}{1-\frac{2}{n^n\left(\frac{1}{n}\right)^{n-1}}}=\frac{n^3}{n-2}.
$$
As an immediate consequence of this and (6), we obtain an upper bound for the best constant $$\lambda_n\le \frac{n^2}{n-2}.\eqno(7)$$

We want to use a well known result
in analysis, which states that a continuous function over a compact set achieves its
minimum (and maximum) values at certain points. For this purpose, let us change function $g(x_1,\ldots,x_n)$, to new function $g_1$ so that $g_1$ is defined also at point $P_0\left(\frac{1}{n},\frac{1}{n},\ldots,\frac{1}{n}\right)$ and points of $\partial C$, and  $g_1$ is continuous in compact set $\overline{C}=C\cup\partial C$:
$$
g_1(x_1,\ldots,x_n)=\begin{cases}

\frac{\pi}{2}, & \text{ if } \prod\limits_{i=1}^n{x_i}= 0;\\
\tan^{-1} {\frac{n^3}{n-2}}, & \text{ if } (x_1,x_2,\ldots,x_n)= \left(\frac{1}{n},\frac{1}{n},\ldots,\frac{1}{n}\right);\\
\tan^{-1} {g(x_1,\ldots,x_n)}, & \text{ otherwise.}
\end{cases}
$$
Since $g_1$ is a continuous function in compact $\overline{C}$, $g_1$ reaches its extreme values somewhere in $\overline{C}$. Obviously, $g_1$ reaches its maximum value $\frac{\pi}{2}$ at the boundary points $\partial C$ where $\prod\limits_{i=1}^n{x_i}=0$, and the minimum value at a point of $C$. The minimum of $g$ is achieved at the same point of $C$ if the minimum point is different from $P_0\left(\frac{1}{n},\frac{1}{n},\ldots,\frac{1}{n}\right)$. In any case, $\inf\limits_{\textbf{x}\in C}{g}=\tan\left(\min\limits_{\textbf{x}\in \overline{C}}{g_1}\right)$. We use an optimization argument similar to \cite{sato, mitev} but with 3 variables, to determine where these points must lie. This method can also be used for other inequalities involving only symmetric polynomials $s_1$, $s_{n-1}$, and $s_n$.

Let $P(x_{1},x_{2},\ldots ,x_{n})$ be a minimum point of $g_1$. Select any 3 of the coordinates of $(x_{1},x_{2},\ldots ,x_{n})$, say $x_1$, $x_2$, and $x_3$. Let us assume that $x_1x_2x_3=\alpha$ and $x_1+x_2+x_3=\beta$. Since $P\in C$, $\alpha,\beta>0$. Also, by AM-GM inequality $\beta^3\ge 27\alpha$ and it is known that if $\beta^3= 27\alpha$ then $x_1=x_2=x_3$. So, suppose that $\beta^3>27\alpha$. Let us now take arbitrary positive numbers  $x,y,z$ such that $xyz=\alpha$ and $x+y+z=\beta$. Without loss of generality we can assume that $x\le y\le z$. Since $x+z=\beta-y$ and $xz=\frac{\alpha}{y}$, the numbers $x$ and $z$ are the solutions of quadratic equation $\delta^2+(y-\beta)\delta+\frac{\alpha}{y}=0$. If we take $y=t$ then we obtain parametrization of the curve obtained by intersection of plane $x+y+z=\beta$ and surface $xyz=\alpha$:
$$
x=\frac{-t+\beta\pm\sqrt{(t-\beta)^2-\frac{4\alpha}{t}}}{2},\ y=t,\ z=\frac{-t+\beta\mp\sqrt{(t-\beta)^2-\frac{4\alpha}{t}}}{2}.
$$
Parameter $t$ changes in interval $[t_1,t_2]$, where $t_1$ and $t_2$ are the zeros of cubic $\kappa(x)=t(t-\beta)^2-4\alpha$ in intervals $\left(0,\frac{\beta}{3}\right)$ and $\left(\frac{\beta}{3},\beta\right)$, respectively. The third zero $t_3$ of $\kappa(x)$ satisfies $t_3>\beta$ and therefore $t_3\not\in[t_1,t_2]$. Since we are interested only with case $x\le y\le z$, we will take one half of this curve  (see Fig. 2)
$$
x=\frac{-t+\beta-\sqrt{(t-\beta)^2-\frac{4\alpha}{t}}}{2},\ y=t,\ z=\frac{-t+\beta+\sqrt{(t-\beta)^2-\frac{4\alpha}{t}}}{2},
$$
and in smaller interval $[t^*_1,t^*_2]$, where $t^*_1$ and $t^*_2$ are the zeros of cubic $\kappa^*(x)=\kappa(x)-t(3t-\beta)^2$ in intervals $\left(t_1,\frac{\beta}{3}\right)$ and $\left(\frac{\beta}{3},t_2\right)$, respectively. Again, since the third zero $t^*_3$ of $\kappa^*(x)$ satisfies $t^*_3>\beta$, $t^*_3\not\in[t_1,t_2]$. Note that if $t=t^*_1$, then $x=y$, and if $t=t^*_2$ then $y=z$. Consider sum $\frac{1}{x}+\frac{1}{y}+\frac{1}{z}$  and note that
$$
\frac{1}{x}+\frac{1}{y}+\frac{1}{z}= \frac{1}{y}+\frac{x+z}{xz}=\frac{1}{y}+\frac{\beta-y}{\alpha/y}=\frac{1}{t}+\frac{\beta t-t^2}{\alpha}.
$$
\begin{figure}[htbp]
\centerline{\includegraphics[scale=.5]{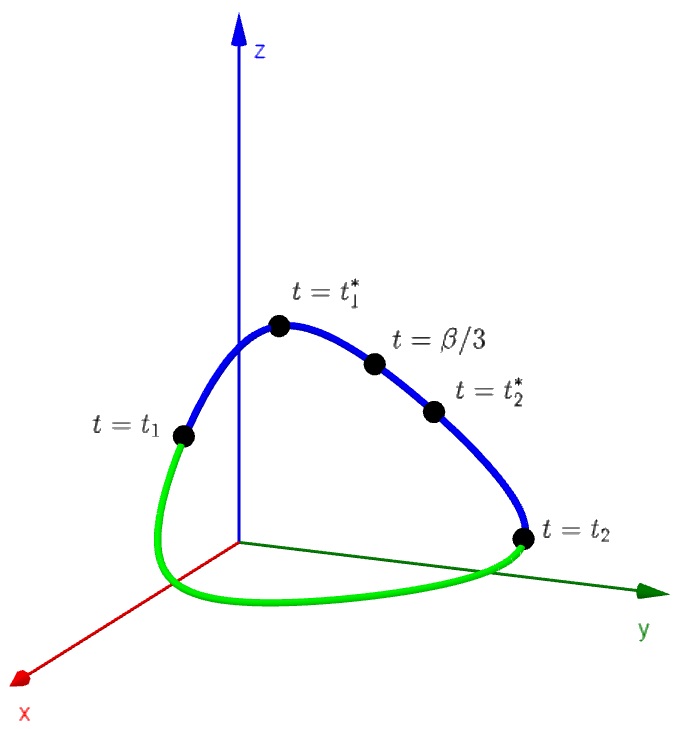}}
\label{fig2}
\caption{Parametric space curve (blue and green) representing intersection of plane $x+y+z=\beta$ (not shown) and surface $xyz=\alpha$ (not shown).}
\end{figure}
Denote $h(t)=\frac{1}{t}+\frac{\beta t-t^2}{\alpha}$, where $t\in[t^*_1,t^*_2]$. So, if $t\in(t^*_1,t^*_2)$, then $h'(t)=-\frac{1}{t^2}+\frac{\beta}{\alpha}-\frac{2t}{\alpha}=\frac{(z-t)(t-x)}{xy^2z}> 0$, and $h'(t^*_1)=h'(t^*_2)=0$. Consequently, $h(t)$ attains its minimum and maximum in interval $[t^*_1,t^*_2]$ at endpoints $t^*_1,\ t^*_2$, respectively.
We are interested in making $h$ smaller, which happens when sum $\frac{1}{x}+\frac{1}{y}+\frac{1}{z}$ is smaller. So, the minimum of $\frac{1}{x_{1}}+\frac{1}{x_{2}}+\frac{1}{x_{3}}$ is reached when $x_{1}=x_{2}\le x_{3}$. Since the coordinates $x_1$, $x_2$, and $x_3$ were chosen arbitrarily, these results hold
for any trio of coordinates. Therefore, the left side of (6) is minimal only when there are at most 2 distinct numbers in the set $\{x_{1},x_{2},\ldots ,x_{n}\}$ Furthermore if the two numbers are distinct then the smaller one is repeated $n-1$ times in $\{x_{1},x_{2},\ldots ,x_{n}\}$ i.e. $x_{1}=x_{2}=\ldots =x_{n-1}\le x_{n+1}$. Consequently, in (6) we can restrict ourselves only to the case where  $x_{1}=x_{2}=\ldots =x_{n-1}=x$, $x_{n}=1-(n-1)x$, where $0<x\le\frac{1}{n}$. By substituting these in (6) and simplifying, we obtain
$$
\frac{n^2}{1-(n-1)\frac{(nx-1)^2}{((n-1)-n(n-2)x)(1-n^nx^{n-1}(1-(n-1)x))}}\geq \lambda.
$$
We will study the part of the denominator which is dependent on $x$, and for simplicity put $t=nx$. So,
$$
\frac{(t-1)^2}{((n-1)-(n-2)t)(1-nt^{n-1}+(n-1)t^n)}=
$$
$$
\frac{1}{((n-1)-(n-2)t)(1+2t+3t^2+\ldots+(n-1)t^{n-2})}.
$$
Denote the polynomial in the denominator by
$$
p_n(t)=((n-1)-(n-2)t)(1+2t+3t^2+\ldots+(n-1)t^{n-2}),
$$
where $0\le t\le 1$. By taking the derivative and simplifying, we obtain
$$
p_n^\prime (t)=(n-1)\left(1\cdot 2+2\cdot 3t +3\cdot 4t^2+\ldots+(n-2)\cdot (n-1)t^{n-3}\right)
$$
$$
-(n-2)\left(1^2+2^2t +3^2t^2+\ldots+(n-1)^2t^{n-2}\right)
$$
$$
=1\cdot n+2\cdot (n+1)t +3\cdot (n+2)t^2+\ldots+(n-2)\cdot (2n-3)t^{n-3}-(n-2)(n-1)^2t^{n-2}.
$$
Since $p_n^\prime (0)=n>2$ and $p_n^\prime (1)=-\frac{n(n-1)(n-2)}{6}<0$, there is at least 1 zero of polynomial $
p_n^\prime (t)
$ in interval $(0,1)$. On the other hand, by Descartes' rule of signs (see p. 247 in \cite{kurosh}, or p. 28 in \cite{prasolov}) the number of positive zeros of $p_n^\prime (t)$ does not exceed the number of sign changes in the sequence of coefficients of $p_n^\prime (t)$, which is 1. So, $p_n^\prime (t)$ has exactly one zero $t_n$ in $[0,1]$, which is also maximum point of $p_n(t)$. This means that there is exactly one point $x=\frac{t_n}{n}$ in $\left(0,\frac{1}{n}\right)$, such that
$x_{1}=x_{2}=\ldots =x_{n-1}=x$, $x_{n}=1-(n-1)x$ makes the left side of (6) minimal. This minimal value is also the best constant for (1):
$$
 \lambda_n=\frac{n^2}{1-\frac{n-1}{p_n(t_n)}}. \eqno(8)
$$
For $n=3,4,5,$ and $6$ it is possible to find the exact values of $t_n$ and corresponding $\lambda_n$.
\begin{itemize}
\item if $n=3$, then $
p_3(t)=(2-t)(1+2t),
$ and $
p_3^\prime (t)=3-4t.
$ Therefore, $t_3=\frac{3}{4}$. By (8), the best constant is $\lambda_3=\frac{3^2}{1-\frac{3-1}{p_3(t_3)}}=25$ (see \cite{aliyev1}).

\item if $n=4$, then $
p_4(t)=(3-2t)(1+2t+3t^2),
$ and $
p_4^\prime (t)=4+10t-18t^2.
$ Therefore, $t_4=\frac{5+\sqrt{97}}{18}$. By (8), the best constant is $\lambda_4=\frac{4^2}{1-\frac{4-1}{p_4(t_4)}}=\frac{582\sqrt{97}-2054}{121}\approx 30.423077$ (see \cite{aliyev2}).

\item if $n=5$, then $
p_5(t)=(4-3t)(1+2t+3t^2+4t^3),
$ and $
p_5^\prime (t)=5+12t+21t^2-48t^3.
$ Using Cardano's formula and Maple, we find that $t_5=\frac{\theta+7+241\theta^{-1}}{48},$ where $\theta=\left(8119+48 \sqrt{22535}\right)^{\frac{1}{3}}$. By (8), the best constant is
$
\lambda_5=\frac{5^2}{1-\frac{5-1}{p_5(t_5)}}\approx 40.090307$, which coincides with the value of $\lambda_5$ conjectured in \cite{aliyev2}.

\item if $n=6$, then $$
p_6(t)=(5-4t)(1+2t+3t^2+4t^3+5t^4),\ 
p_6^\prime (t)=6+14t+24t^2+36t^3-100t^4.
$$
Using Ferrari's method and Maple, we find that
$$t_6=\frac{9+\phi+\sqrt{50 \psi+962 -11300 \psi^{-1}+47258 \phi^{-1}}}{100},$$
where
$$
\phi=\sqrt{-50 \psi+481 +11300\psi^{-1}}, \psi=\left(1473+\sqrt{13712905}\right)^{\frac{1}{3}}.
$$
By (8), the best constant is
$
\lambda_6=\frac{6^2}{1-\frac{6-1}{p_6(t_6)}}\approx 52.358913$.
\end{itemize}
For larger values of $n$, we can give some bounds for $\lambda_n$. We already found an upper bound (7). We will now focus on a similar lower bound.

By AM-GM inequality,
$$
\sum\limits_{i = 1}^n{\frac{1}{x_{i}}}\geq \frac{n
}{G_n}, \eqno(9)
$$
where $x_{1},x_{2},\ldots ,x_{n}>0; \sum\limits_{i =
1}^n{x_{i}}=1$,  $G_n=\sqrt[n]{\prod\limits_{i=1}^n{x_{i}}}$, and $n\geq 2$. Let us show that if $\lambda=\frac{n^3}{n-1}$ then
$$\frac{n
}{G_n}\ge\frac{\lambda
}{1+n^{n-2}(\lambda -n^{2})G_n^n},\eqno(10)
$$
where $(x_1,x_2,\ldots,x_n)\ne \left(\frac{1}{n},\frac{1}{n},\ldots,\frac{1}{n}\right)$ and therefore $G_n=\sqrt[n]{\prod\limits_{i=1}^n{x_{i}}}<\frac{1}{n}$. Indeed, we can simplify (10) to
$$
\frac{n^2(1-s^n)}{s(1-s^{n-1})}\ge\frac{n^3}{n-1},
$$
where $s=nG_n<1$. It is easily proved, as we can write it in the following form
$$
1+\frac{1}{s+s^2+\ldots+s^{n-1}}\ge \frac{n}{n-1},
$$
where noting $s<1$ completes proof. From (1), (9), and (10) it follows that if $\lambda\le\frac{n^3}{n-1}$, then (1) holds true. This means that we have now a lower bound for the best constant:
 $$\lambda_n\ge\frac{n^3}{n-1}\eqno(11).$$
Combining (7) and (11) we obtain the following symmetric double inequality.

\textbf{Theorem 1.} If $n>2$, then
 $$\frac{n^3}{n-1}\le\lambda_n\le\frac{n^3}{n-2}\eqno(12).$$

It is possible to improve these estimates in exchange for less elegant formula. For example, if we put $x_{1}=x_{2}=\ldots =x_{n-1}=\frac{1}{n+1}$, $x_{n}=\frac{2}{n+1}$, then we obtain from (6) a new upper bound for the best constant:
$$
\lambda_n\le (n+1)^2\cdot\frac{\frac{1}{2}-\frac{n^n}{(n+1)^n}}{\frac{n+1}{2n-1}-\frac{n^{n-2}}{(n+1)^{n-2}}}.\eqno(13)
$$
One can check that (13) is sharper than (6) for all $n>3$. We can also prove that $\frac{n}{n+1}\le t_n$ or equivalently, $p^{\prime}_n\left(\frac{n}{n+1}\right)\ge0$ for all $n\ge3$. Indeed, $$p^{\prime}_n\left(\frac{n}{n+1}\right)=3 n \left(n+1\right)^{2} \left(\frac{n}{n+1}\right)^{n} \left(\left(1+\frac{1}{n}\right)^{n}-\frac{8}{3}+\frac{1}{n}-\frac{1}{3 n^{2}}\right).$$
For $n=3,4,\ldots,25$ one can check directly that $\left(1+\frac{1}{n}\right)^{n}-\frac{8}{3}+\frac{1}{n}-\frac{1}{3 n^{2}}\ge 0$. For $n>25$, one can use the fact that $\left(1+\frac{1}{n}\right)^{n}>\frac{8}{3}$ and $\frac{1}{n}>\frac{1}{3 n^{2}}$.

\begin{section}
{Appendix}

We will give a proof of (5) here. We can use the optimization argument given after Lemma 1 of the current paper, to maximize the left hand side of (5), while keeping the right hand side of (5) fixed. This is achieved when for any 3 of the coordinates, say $x_1$, $x_2$, and $x_3$, of $(x_{1},x_{2},\ldots ,x_{n})$, $x_1\le x_2=x_3$. So, we can restrict ourselves only to the case where  $x_{1}=x$ and $x_{2}=\ldots =x_{n-1}=x_{n}=\frac{1-x}{n-1}$, where $0<x\le\frac{1}{n}$. For this particular case (5) is transformed to
$$
\frac{1}{x}+\frac{(n-1)^2}{1-x}\leq \frac{(n-1)^{n-1}
}{n^{n-2}x(1-x)^{n-1}},
$$
which can be simplified to correct inequality
$$
(nx-1)^2\left((n-2)(nx)^{n-3}+(n-3)(nx)^{n-4}+\ldots+1\right)\ge 0.
$$
The equality case is possible only when $nx=1$.

Inequality (5) can also be written as homogeneous inequaity $A_n^{n-1}H_n\ge G_n^n$, where $A_n,\ H_n,$ and $G_n$ are, respectively, arithmetic, harmonic, and geometric means of arbitrary positive numbers $x_1,\ldots,x_n$:
$$
A_n=\frac{\sum\limits_{i =
1}^n{x_{i}}}{n},\ H_n=\frac{n}{\sum\limits_{i =
1}^n{\frac{1}{x_{i}}}},\ G_n=\sqrt[n]{\prod\limits_{i =
1}^n{x_{i}}}.
$$
Since $A_n\ge G_n$, automatically $A_n^{l}H_n\ge G_n^{l+1}$ for any real number $l\ge n-1$. It would be natural to ask whether general inequality $A_n^{l}H_n\ge G_n^{l+1}$ can hold true also for some real number $l<n-1$. The answer to this question is negative. A counter example is found if one takes $x_{1}=x_{2}=\ldots =x_{n-1}=1$ and $x_{n}=x$, where $x\rightarrow 0^+$. Indeed,  if $l<n-1$, then $A_n^{l}H_n=O(x)$ and $G_n^{l+1}=O(x^{\frac{l+1}{n}})\gg O(x)$.

\end{section}

\end{document}